\documentclass[12pt]{article}

\usepackage{fullpage}
\usepackage{epsfig}
\usepackage{color}
\usepackage{bbm}

\newtheorem{theorem}{Theorem}    
\newtheorem{lemma}{Lemma}      
\newtheorem{remark}{Remark}

\def\eod{\vrule height 6pt width 5pt depth 0pt}
\newenvironment{proof}{\noindent {\bf Proof:} \hspace{.2em}}
                      {\hspace*{\fill}{\eod}}
                      {}
\newcommand{\os}{\overline{s}}
\newcommand{\oa}{\overline{a}}
\newcommand{\ob}{\overline{b}}
\newcommand{\cT}{{\cal T}}
\newcommand{\Pk}{\mathsf{Park}}

\newcommand{\sT}{\mathsf{sign}(T)}
\newcommand{\hgt}{{\sf ht}}

\newcommand{\comment}[1]{}
\newcommand{\cF}{\cal F}

\begin{document}

\title{Spanning trees in complete uniform hypergraphs and
a connection to $r$-extended Shi hyperplane arrangements}
\author{Sivaramakrishnan Sivasubramanian\\
Institute of Computer Science\\
Christian-Albrechts-University\\
Kiel 24118, Germany\\
email: ssi@informatik.uni-kiel.de
}

\date{\today}

\maketitle

\begin{abstract}
We give a Cayley type formula to count the number of spanning trees in the 
complete $r$-uniform hypergraph for all $r \geq 3$.  
Similar to the bijection between spanning trees of the complete 
graph on $(n+1)$ vertices and {\sl Parking functions} of length $n$, we derive a 
bijection from spanning trees of the complete
$(r+1)$-uniform hypergraph which arise from a fixed $r$-perfect matching (see Section 
\ref{sec:unif}) and {\sl $r$-Parking functions} of an appropriate length.  We 
observe a simple consequence of this bijection
in terms of the number of regions of the $r$-extended Shi hyperplane arrangement 
in $n$ dimensions, $S_n^r$.
\end{abstract}

\section{Introduction}

We give a formula to count the number of spanning trees in the complete 
$r$-uniform hypergraph for all $r \geq 3$ (we call them $r$-spanning trees).  We first 
present the case when $r=3$ where 
we use the Pfaffian Matrix Tree Theorem of Masbaum and Vaintrob \cite{MV-02}.  
Using ideas from that proof, we present our result for $r \geq 4$.

\begin{figure}[h]
\centerline{\psfig{file=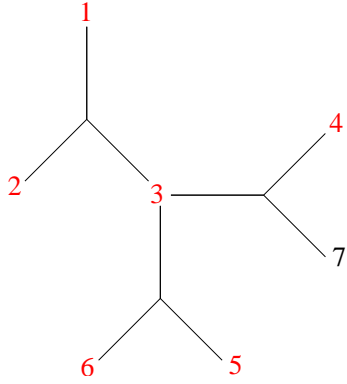,height=5cm}}
\caption{A 3-uniform spanning tree on 7 vertices}
\label{sptree:example}
\end{figure}

For a positive integer $n$, let $[n] = \{1,2,\ldots,n \}$. 
3-spanning trees are defined in the following way.  Let $H=(V,\cF)$ be
a 3-uniform hypergraph.  Consider the
bipartite graph with $V$ on one side,  and $\cF$ on the other.  Each
hyperedge $e=\{a,b,c\}$ is connected to precisely the vertices $a,b$ and $c$.  If the
bipartite graph on some $r$ hyperedges $E = \{e_1,e_2,\ldots,e_r\}$ and all the vertices $V$
is a tree in the graph theoretic sense, then we call the hypergraph $H=(V,E)$ as a 3-spanning tree.  
Figure \ref{sptree:example} shows a 3-spanning tree $T=([7],\{123, 347, 356 \})$ on 7 vertices
where $123$ is an abbreviation for the 3-hyperedge $\{1,2,3\}$ and so on.
We prove the following theorem about 3-spanning trees.

\begin{theorem} 
\label{thm:count}
The number of 3-spanning trees of the complete 3-uniform hypergraph on 
$n=2k+1$ vertices is $1*3\cdots*(2k-1)*(2k+1)^{k-1}$.
\end{theorem}

By counting the number of edges in the bipartite graph representation of a 3-spanning
tree in two ways, it is clear that for spanning trees to exist in a 3-uniform 
hypergraph on $n$ vertices, $n$ has to be odd.  The same
double counting argument shows that the number of hyperedges in any 3-spanning tree 
on $n=2k+1$ vertices is $k$.  Thus  Theorem \ref{thm:count} is similar
to Cayley's Theorem for counting spanning trees of the complete graph on $n$ vertices.
If we interpret the term occuring in
Cayley's Theorem as $v^{e-1}$ with $v$ being the number of vertices spanned
and $e$ being the number of edges in any spanning tree, we see that Theorem \ref{thm:count} has
a similar term $(2k+1)^{k-1}$.   There is an additional multiplicative
term of $1*3*\cdots*(2k-1)$.  This term is the number of perfect matchings in the
complete graph on $n-1 = 2k$ vertices and from the 
proof of Theorem \ref{thm:count}, one sees why this term arises.

Let $r \geq 4$ be a positive interger.  $r$-spanning trees in $r$-uniform hypergraphs are defined 
analogously.  Figure \ref{fig:4sptree} shows a 4-spanning tree.  Henceforth, when we talk 
about $r$-spanning trees, we omit mentioning that the underlying graph is 
the complete $r$-uniform hypergraph.  
We prove a similar result about the number of $r$-spanning trees on $n = (r-1)k + 1$ 
vertices.  For $r \geq 3$, it can be checked that the number of $(r-1)$-perfect matchings 
of the complete $(r-1)$-uniform hypergraph on $(r-1)k$ vertices is
$rPM = {{(r-1)\cdot k -1} \choose r-2} * {{(r-1)\cdot (k-1) -1} \choose r-2} * \cdots * 
{r-2 \choose r-2}$.
The proofs of both Theorems \ref{thm:count} and \ref{thm:rcount} are in Section \ref{sec:unif}.

\begin{theorem}
\label{thm:rcount}
For $r \geq 3$, the number of $r$-spanning trees on $n=(r-1)k+1$ vertices 
is $rPM * n^{k -1}$.
\end{theorem}

In Section \ref{sec:egf}, we give an exponential generating function for the number of 
rooted $r$-spanning trees on $[n]$ .  These are analogs of the famous
relation $D(x) = x \exp(D(x))$ where $D(x)$ is the exponential generating function 
for rooted spanning trees.  Though one can derive Theorems
\ref{thm:count} and 
\ref{thm:rcount} from this (ie the total count), we present a non generating function 
proof of both these theorems in Section \ref{sec:unif} as we get extra information about the number 
of $r$-spanning trees
on $[n]$, when we fix an $(r-1)$-perfect matching on $n-1$ vertices.  
This count of a subset of the set of all
$r$-spanning trees will be used in Sections \ref{sec:park} and \ref{sec:shi}.

In Section \ref{sec:park}, we give a bijection between $r$-parking functions of length $k$
and $(r+1)$-spanning trees on $rk+1$ vertices, which arise from a fixed 
$r$-perfect matching on $rk$ vertices.
We recall the definition of {\sl parking functions} of length $k$.  
There are $k$ cars $1,2,\ldots,k$ and $k$ parking spaces $0,1,\ldots,k-1$, in this order.
All cars enter the street at the end close to the zeroth parking slot.  They enter
in increasing order and each car $i$ has its preferred parking slot $a_i$.  Car $i$ drives to
slot $a_i$ and parks there if that slot is free.  If not, it tries the next higher parking slot 
and so on till it gets an empty parking slot.  A sequence 
$\oa = (a_1,a_2,\ldots,a_k) \in  \mathbbm{Z}_{\geq 0}^k$ is said 
to be a 
{\sl parking function} of length $k$ if all cars are able to park following the above rules.  
Spanning trees of the complete graph on 
$k+1$ vertices are in a bijection with {\sl parking functions} of length $k$ (see \cite{EC2}).
There is an alternative characterization of parking functions.  Let 
$\oa = (a_1,a_2,\ldots,a_k) \in \mathbbm{Z}_{\geq 0}^k$ 
and let $\ob = (b_1,b_2,\ldots,b_k)$ 
where $b_1 \leq b_2 \leq \cdots \leq b_k$ be the weakly increasing
rearrangement of $\oa$. 
Then $\oa$ is a parking function iff $b_i \leq (i-1)$ for all $i$.  This
algebraic definition of Parking functions was generalised in the following manner to yield
$r$-parking functions for a positive integer $r$ (see \cite{stan-rota}).  A sequence 
$\oa \in \mathbbm{Z}_{\geq 0}^k$
is called an $r$-parking function iff its weakly increasing rearrangement $\ob$
satisfies $b_i \leq r(i-1)$ for all $i$.  
For $r,k \geq 1$, let $\Pk_k^r$ be the set of $r$-parking functions of length $k$.
We prove the following theorem.

\begin{theorem}
\label{thm:bij}
For all $r \geq 1$, there is a bijection between the set of $(r+1)$-spanning trees on $n=rk+1$ 
vertices which arise from a fixed $r$-perfect matching and the set $\Pk_k^r$.
\end{theorem}

Our bijection is very similar to that of Chebekin and Pylyavskyy \cite{CP} and uses the 
{\sl breadth-first search} order of vertices in a rooted $r$-spanning tree.

In Section \ref{sec:shi}, we point out a simple connection to the number of regions of the 
$r$-extended Shi hyperplane arrangement in $m$ dimensions, denoted $S_m^r$.  This connection is a 
consequence of Theorem \ref{thm:bij} and
Theorem 2.1 in \cite{stan-rota}.  We recall the definition of $S_m^r$.  
It is given by the following set of hyperplanes in $\mathbbm{R}^m$ 
$$ x_i - x_j = -r+1,-r+2,\ldots,r, \mbox{ for } 1 \leq i < j \leq m$$

It has ${m \choose 2} * 2r$ hyperplanes.  The number of regions of the 1-extended Shi 
hyperplane arrangement in $n$ dimensions is identical
to the number of spanning trees of the complete graph on $n+1$ vertices (see \cite{stan-rota}).  
We prove an analog of this for higher values of $r$.

\begin{theorem}
\label{thm:easy}
Let $r,k \geq 1$.  The number of regions of the $r$-extended Shi arrangement $S_{k}^r$ is equal 
to the number of $(r+1)$-spanning trees on $n=rk+1$ vertices arising 
from a fixed $r$-perfect matching on $n-1$ vertices.
\end{theorem}

\section{Counting $r$-spanning trees}
\label{sec:unif}

We first prove our result for 3-spanning trees and then for $r$-spanning trees where
$r \geq 4$.  

\subsection{Counting 3-spanning trees}


We review briefly the theorem of Masbaum and Vaintrob
\cite{MV-02} where they enumerate with a $\pm$ sign all 3-spanning trees of a 3-uniform
hypergraph.  We only need the case when the 3-uniform hypergraph is complete.
As remarked, for 3-spanning trees to exist in 3-uniform 
hypergraphs, the number of vertices has to be odd.  Let $n=2k+1$ 
and let $ST_n$ be the set of 3-spanning trees on $n$ vertices.
For $i,j,k \in [n]$, all three indices being distinct, let 
$x_{i,j,k}$ be a variable with the following `sign' property.  For all 
$\sigma \in S_3$, $x_{\sigma(i),\sigma(j),\sigma(k)} = \mbox{sign}(\sigma) 
x_{i,j,k}$.

Associate an $n \times n$ matrix $A = (a_{i,j})$ with the complete 3-uniform 
hypergraph where
$a_{i,j} = \sum_{k=1;k\not=i,j}^n x_{i,j,k}$.  Due to the sign property, this
matrix is skew-symmetric.  Let $A_i$ be the submatrix of $A$ obtained by omitting the
$i$-th row and the $i$-th column.  

The theorem of Masbaum and Vaintrob says that the pfaffian of $A_i$ gives 
a $\pm 1$ signed enumeration of all 3-spanning trees $T \in ST_n$.
The coefficient $\pm1$ arises dues to two multiplicative factors: a $\pm 1$ sign 
denoted $\sT$ which is the sign of a permutation $\pi$ of the vertices 
obtained by embedding $T$ in the plane; 
and a product of $x_{i,j,k}$'s for each hyperedge $\{i,j,k\}$ of $T$.
Each hyperedge variable $x_{i,j,k}$ has to be ordered acording to $\pi$.
Due to the `sign property', each such hyperedge might get a $\pm 1$.  Hence 
multiplication yields an overall $\pm 1$ coefficient.  One thus gets a term $y_T$ with a 
$\pm 1$ coefficient for each 3-spanning tree $T$.  
See Reiner and Hirschman \cite{RH-02} for an exact procedure to obtain $y_T$ from $T$.

Let $PM_{2k}$ be the set of perfect matchings in the complete graph on $2k$ vertices.
Clearly, $|PM_{2k}| = 1*3*\cdots*(2k-1)$.  It is known that the pfaffian of a skew 
symmetric $2k \times 2k$ matrix $A$ of is given by the expression (see \cite{RH-02})

\begin{equation}
\label{eqn:pfaff}
Pf(A) = \sum_{M \in PM_{2k}} (-1)^{\mathrm{cross}(M)} \prod_{i,j \in M, i<j }a_{i,j}
\end{equation}

where $\mbox{cross}(M)$ is the number of edges which cross the perfect 
matching $M$.  Formally, $\mbox{cross}(M) = |\{i<j<k<l: \{i,k\}, \{j,l\} \in M  \}|$.

In the above expression, we call the terms arising from a fixed perfect matching $M$ as
the ``terms of $M$ in the pfaffian expansion".  The Pfaffian Matrix Tree Theorem is:

\begin{theorem}[Masbaum and Vaintrob \cite{MV-02}]
Let $A$ be the matrix defined above and let $n$ be odd.  For all $1 \leq i \leq n$,
$$ Pf(A_i) = (-1)^{i-1} \sum_{T \in ST_n} y_T
$$
\end{theorem}

\begin{proof}(Of Theorem \ref{thm:count})  

Consider the matrix $A_n$ obtained by deleting the $n$-th row and the $n$-th column of $A$
and expand the pfaffian of $A_n$ as in equation \ref{eqn:pfaff}.  We have terms for each perfect
matching $M \in PM$.  
It is conceivable that some 3-spanning trees $T$ may occur as terms of many perfect 
matchings, and they sum up nicely to yield a $\pm 1$ coefficient for each such $T$.  We show 
below that
this cannot happen.  This observation is what is used to generalise to get the result
for the complete $r$-uniform hypergraph for $r \geq 4$.

We show that each 3-spanning tree appears exactly once in the pfaffian expansion of 
equation \ref{eqn:pfaff} in the term corresponding to some perfect matching.  We recall
that $ST_n$ is the
set of spanning trees of the complete 3-uniform hypergraph on $n$ vertices, 
and $PM_{2k}$ is the set of perfect matchings of the complete graph on $2k$ vertices.

\begin{lemma}  
\label{lemma:perfect}
There is a one-to-one mapping $f: ST_n \mapsto PM_{2k}$ such that for all $T \in ST_n$,
the term $y_T$ appears in the pfaffian expansion of $A_n$ only among the terms of the 
perfect matching $f(T)$.
\end{lemma}

\begin{proof}  We give an algorithm to obtain the mapping $f$.
Let $T \in ST_n$ with $T=(V,E)$.  
Since $T$ spans all the vertices, vertex $n$ has to appear in at 
least one hyperedge.  As each hyperedge 
has exactly 3 vertices, deleting the vertex $n$ from all hyperedges of $E$ will 
result in some hyperedges having size 2 and others retaining their size. (At least one 
hyperedge will have its size reduced by this procedure.)  All the currently size 2 hyperedges 
are put in $f(T)$.  The vertices newly added to $f(T)$ are termed ``matched".  If $f(T)$ is a perfect 
matching, then we are done.  Otherwise, there still remain some hyperedges in $E$ which have 
size 3.  Delete all the {\sl matched} vertices.  This will again reduce the size of some 
hyperedges to 2.  These size 2 hyperedges are added to $f(T)$ and we iterate.  It is clear that 
since we obtain at least one matching edge in each iteration, that in at most $k$ iterations, 
we will terminate with a perfect matching on the vertex set $[2k]$.  
It is also clear that we end up with exactly one perfect matching of $[2k]$.  The perfect 
matching we construct has the property that from every hyperedge of $T$, we have 2 
vertices matched.  For example, on the spanning tree of Figure \ref{sptree:example}, the
algorithm will output the perfect matching $\left\{ \{3,4\}, \{1,2\}, \{5,6\} \right\}$.

We need to check that $y_T$ occurs
in the pfaffian expansion corresponding to the perfect matching $f(T)$.  This is easy as
(ignoring the $\pm$ sign of $y_T$,) all hyperedges of $T$ have exactly two vertices in the perfect 
matching $f(T)$. We also need to check that
$y_T$ does not occur among terms of any other perfect matching in the pfaffian expansion of $A_n$.  
To see this,
we note that two vertices from all hyperedges in $E$ need to occur in any perfect matching $M'$ for
$y_T$ to occur among the terms of $M'$ and that this happens only for $f(T)$.

\end{proof}

\begin{remark}  
\label{rem:pm}
It is easy to see from the above algorithm that each spanning tree $T$ has
at least one leaf (ie a hyperedge with two vertices of degree 1) which is matched in the 
perfect matching.  This will be used in the bijection of Lemma \ref{lem:prufer}.  Suppose we call 
the edge of the perfect matching
$e = \{a,b\}$, then we note that there is a unique hyperedge $f \in E(T)$ such that $f=\{a,b,\ell \}$
for some $\ell \in ([n] - \{a,b\})$.
\end{remark}

\begin{remark}
\label{rem:endpt}
From the above algorithm, it is clear that each edge $e=\{a,b \}$ of the perfect matching
$f(T)$ is in a unique hyperedge $\{a,b,x \}$ of $T$.
\end{remark}

Since each 3-spanning tree $T$ comes up in the pfaffian expansion in exactly one perfect
matching, we infer that each
3-spanning tree occurs exactly once among the terms in the pfaffian expansion of $A_n$.  
Hence, only terms
corresponding to non spanning trees (ie terms with cycles in it) get cancelled and
Reiner and Hirschman \cite{RH-02} exhibited a sign reversing involution cancelling 
exactly these cyclic terms.
Lemma \ref{lemma:perfect} shows that each 3-spanning tree occurs in the 
pfaffian expansion in exactly one perfect
matching.  Below we count the number of 3-spanning trees $t_M$ which get mapped to a fixed
perfect matching $M$ under $f$.   Theorem \ref{thm:count} follows by adding $t_M$ over all
perfect matchings $M$ in $PM_{2k}$.

\begin{lemma}  
\label{lem:prufer}
Let $n = 2k+1$ and let $M \in PM_{2k}$.  Under $f$, $M$ gives rise to $(2k+1)^{k-1}$ 3-spanning 
trees on the vertex set $[n]$.
\end{lemma}

\begin{proof}
We give a {\sl Prufer} type bijection.  We prove the lemma for the perfect matching 
$M = \{ \{1,2\},\{3,4\},\ldots,\{2k-1,2k\} \}$.  It will be clear that a similar proof works
for other perfect matchings.  We give a bijection between 3-spanning trees which arise 
from the perfect matching $M$ to the set of all strings 
$\os = (s_1,s_2,\ldots,s_{k-1})$ with $k-1$ coordinates where for all $1 \leq i \leq
k-1$, $1 \leq s_i \leq n$.

For one direction, given a 3-spanning tree $T$, by Remark \ref{rem:pm}, we know that there is at least 
one {\sl leaf hyperedge} (ie a hyperedge which contains an edge $e$ of $M$ as degree one
vertices; we refer to the edge $e \in M$ as a {\sl leaf edge}).  
If there are many such leaf hyperedges, we choose a total order $\pi$ on the {\sl leaf edges}
and pick the leaf hyperedge with the smallest (wrt $\pi$) leaf edge.  For the remaining 
part of the proof, we use the total order $\{1,2\} <_{\pi} \{3,4\} <_{\pi} \cdots <_{\pi}
\{2k-1,2k\}$ on the edges of $M$.  We pick one leaf hyperedge in every iteration.
Put $i=1$, let $\ell_i$ be the picked leaf hyperedge and let its lead edge be
$e_i = \{a_i,b_i\}$.  
Let $s_i$ be connection point of $\ell_i$ (ie, let $\ell_i = \{s_i,a_i,b_i\}$ is the unique 
hyperedge containing both $a_i$ and $b_i$).  
Delete both the vertices $a_i$ and $b_i$, increase $i$ by 1 and repeat.  Note that there 
will again exist at least one leaf hyperedge in the deleted subhypergraph.  This gives us a 
sequence of $k-1$ numbers, each in the range 1 to $n$.  For example, when we start with the 
tree in Figure \ref{fig:prufer-one-direc}, we get the sequence $\os=(3,3,4)$.

\begin{figure}[h]
\centerline{\psfig{file=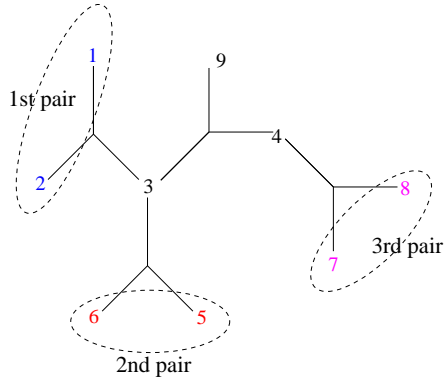,height=5cm}}
\caption{One direction of the Prufer type correspondence.}
\label{fig:prufer-one-direc}
\end{figure}

For the converse, we group the vertices of the graph into components based on 
the perfect matching (see Figure \ref{fig:components}).  Thus, we begin with $k+1$ components 
$\{1,2\}, \{3,4\}, \ldots, \{2k-1,2k\},n$.  As shown in Figure \ref{fig:components}, we need to 
assign values to the variables $x_1,x_2, \ldots,x_k$ such that the resulting 3-uniform hypergraph
is a 3-spanning tree.  
If $1 \leq a < n$ is a number, then let $m(a)$ be the edge of the perfect matching 
$M$ which contains $a$ (ie $m(5) = m(6) = \{5,6\}$ and so on).  If $a = n$, then $m(a)$ is 
undefined. 

\begin{figure}[h]
\centerline{\psfig{file=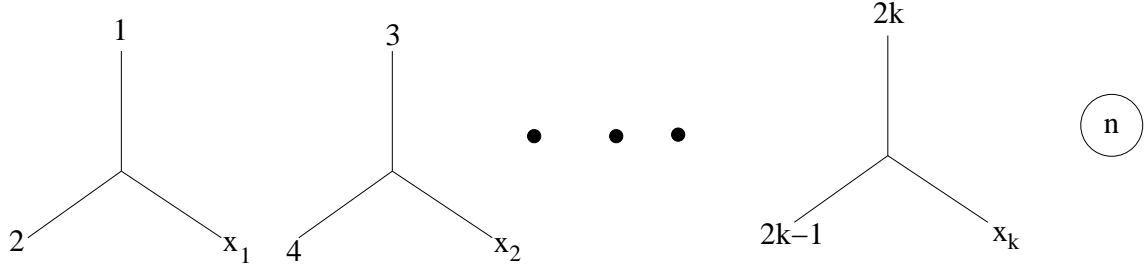,height=3.5cm}}
\caption{The components of $M$ at the beginning.}
\label{fig:components}
\end{figure}
 
Let $\os = (s_1,s_2,\ldots,s_{k-1})$ be the given sequence of numbers.  We add 
hyperedges sequentially.  Initially, all the components are marked ``unfinished".  Starting 
from $i=1$, let $b_i$ be a smallest edge of the perfect matching $M$ which does not occur as 
$m(s_j)$ for $j \geq i$.  We mark $b_i$ as ``finished", add the hyperedge $\{b_i,s_i \}$,
increase $i$ by 1 and iterate.  It is easy to check that we end up with two unfinished 
components.  
We add them as a hyperedge.  For example, the above procedure for $n=9$, on the 
sequence  $\os =(3,3,4)$
will yield the spanning tree of Figure \ref{fig:prufer-one-direc}.
\end{proof}

The above proof clearly works for any perfect matching (we also need a total order between 
the edges of the perfect matching).  Summing over all perfect
matchings in $PM_{2k}$ completes the proof of Theorem \ref{thm:count}.
\end{proof}

\begin{remark}
A similar theorem is true for any 3-uniform hypergraph $H$, though the number of 3-spanning
trees arising from a perfect matching $M$ may depend on $M$ and $H$.
\end{remark}

\subsection{Counting $r$-spanning trees } 
\label{subsec:higher}
In this subsection, we consider $r$-spanning trees where $r \geq 4$.  It can be 
checked that $r$-spanning trees on $n$ vertices exist only when $n \equiv 1 (mod \>\> r-1)$.

\begin{figure}[h]
\centerline{\psfig{file=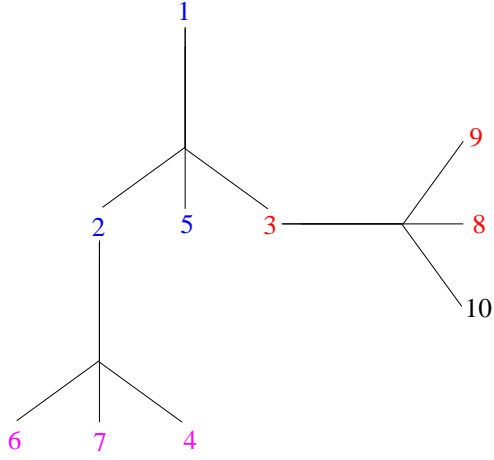,height=6cm}}
\caption{A spanning tree of the 4-uniform complete graph on 10 vertices.}
\label{fig:4sptree}
\end{figure}

Let $T$ be an $r$-spanning tree on $n=(r-1)k+1$ vertices.  
As before, the vertex $n$ lies in at least one hyperedge.  Using 
the same deletion process of Lemma \ref{lemma:perfect},
we see that $T$ arises from exactly one $(r-1)$-perfect matching on the vertex
set $[(r-1)k]$.  We illustrate this on the 4-spanning tree $T$ shown in 
Figure \ref{fig:4sptree}.  The 3d 
matching we get from $T$ is $\{3,8,9\}, \{1,2,5\}$ and $\{4,6,7\}$.  Clearly each $r$-spanning
tree $T$ gives rise to one fixed $(r-1)$-perfect matching and we need to count the number of 
$r$-spanning trees which arise from a fixed $(r-1)$-perfect matching.  

\begin{proof} (Of Theorem \ref{thm:rcount})
A {\sl Prufer} type bijection works in this case as well to show that for a 
fixed $(r-1)$-perfect matching $M$ on $[(r-1)k]$, the number of $r$-spanning trees 
that arise from it equals $(rk+1)^{k-1}$.  
\end{proof}

\section{Exponential generating functions}
\label{sec:egf}

In this section, we present an exponential generating function for the number
of rooted $r$-spanning trees.  We begin with the case $r=3$.  Let $T_n$ be the set of 
rooted
3-spanning trees on $[n]$ and let $t_n$ be the number of rooted 3-spanning trees.  
We set $t_0 = 0; t_1 = 0$.  We note that 
$t_n = 0$ for all even $n$.  Let $$ T(x) =  \sum_{n \geq 0} t_n \frac{x^n}{n!}$$
be the exponential generating function for the sequence $t_n$.  

Let $pm_k$ be the number of perfect matchings of the complete graph on $[k]$.
Let $$E_{pm}(x) = \sum_{n \geq 0} pm_{n} \frac{x^n}{n!} $$
be the exponential generating function of the sequence $pm_k$.  We prove the following
theorem.

\begin{theorem}
\label{thm:egf}
The exponential generating functions $E_{pm}(x)$ and $T(x)$ satisfy 
\begin{equation}
\label{eqn:3egf}
T(x) = x E_{pm}(T(x))
\end{equation}
\end{theorem}

\begin{proof}
We give a 5-step procedure to build rooted 3-spanning trees on $[n]$.  It is easy to
see that Theorem \ref{thm:egf} is equivalent to this procedure. To build a rooted 
3-spanning tree $T$ on $[n]$, we proceed as follows:

\begin{enumerate}
\item Pick $rt$ in $T$ to serve as its root.
\item Partition $[n] - \{rt\}$ into an even number of non empty blocks $B_1,B_2,\ldots,B_{2m}$.
\item Choose a rooted 3-spanning tree $T_i$ on each of the blocks $B_i$ for $1 \leq i \leq 2m$.
\item Choose a perfect matching $\{a_i, b_i\}$ for $i=1,2, \ldots, m$ of the roots of
  $T_1,T_2, \ldots, T_{2m}$.
\item Add $m$ 3-hyperedges $\{rt,a_i,b_i\}$ for $i=1,\ldots,m$.
\end{enumerate}

\end{proof}

One can infer Theorem \ref{thm:count} from the above by applying Lagrange Inversion 
Formula (see \cite{EC2}).  It is easy to see from Equation (\ref{eqn:3egf}) that $t_n = 0$ for
even values of $n$.  Below, we compute $t_n$ for odd $n$. 
 
$$
 	[x^{2k+1}]T(x)    =   \frac{1}{2k+1} [y^{2k}](E_{pm}(y))^{2k+1} 
$$
 
 One can check that 
 $$ [y^{2k}](E_{pm}(y))^{2k+1} = \left( \frac{2k+1}{2} \right)^k \frac{1}{k!} $$
 
Thus, $\frac{t_{2k+1}}{(2k+1)!} = \frac{1}{(2k+1) k!} \left( \frac{2k+1}{2} \right)^k $, or 
$t_{2k+1} = \frac{(2k)!}{k! 2^k} (2k+1)^k = 1 \cdot 3 \cdots (2k-1) (2k+1)^k$.  But
$t_{2k+1}$ is the number of rooted 3-spanning trees on $[2k+1]$ which is $(2k+1)$ times
number of 3-spanning trees.  That completes another proof of Theorem \ref{thm:count}. 

\subsection{Exponential generating function for rooted $r$-spanning trees}

For counting $r$-spanning trees, one modifies the Steps 2 onwards of the above procedure.
To construct an $r$-spanning tree $T$ on $n=(r-1)k + 1$ vertices, we use the procedure below.

\begin{enumerate}
\item Pick $rt$ in $T$ to serve as its root.
\item Partition $[n] - \{rt\}$ into $m$ non empty blocks $B_1,B_2,\ldots,B_m$, where 
$m = p(r-1)$ for a positive integer $p$.
\item Choose a rooted $r$-spanning tree $T_i$ on each of the blocks $B_i$ for $1 \leq i \leq m$.
\item Group the $p(r-1)$ roots of $T_1, T_2, \ldots, T_m$ into $p$ blocks each of size $r-1$.
ie choose an $(r-1)$-perfect matching $R = \{r_1, r_2, \ldots, r_p \}$ of size $p$ where
$r_i = \{ a_1^i, a_2^i,\ldots,a_{r-1}^i \}$ of the roots of $T_1,T_2, \ldots, T_m$.
\item Add $p$ $r$-hyperedges $\{rt\} \cup r_i$ for $i=1,\ldots,p$.
\end{enumerate}

For $r \geq 3$, let $r$-$pm_n$ be the number of $r$-perfect matchings on $[n]$.  Let 
$$ 
E^r_{pm}(x) = \sum_{n \geq 0} \mbox{$r$-$pm_n$} \frac{x^n}{n!}
$$ 
be the exponential generating function of the sequence $r$-$pm_n$.  Let $t^{r+1}_n$ be the 
number of rooted $(r+1)$-spanning trees on $[n]$ and let 
$$ T^{r+1}(x) = \sum_{n \geq 0} t^{r+1}_n \frac{x^n}{n!} $$
be the exponential generating function of the sequence $t^{r+1}_n$.  It is easy to see
from the above procedure that the following.
 
\begin{theorem} For all $r \geq 2$, 
 $$ T^{r+1}(x) =  x E^r_{pm}( T^{r+1}(x) )  $$
\end{theorem}

\section{$r$-Parking functions}
\label{sec:park}
In this section, we prove Theorem \ref{thm:bij}.  We prove the theorem for 
$r=2$.  The proof is identical for higher values of $r$.

\subsection{Connection to 3-spanning trees}

Our proof of Theorem \ref{thm:bij} closely mimicks that of Chebekin and 
Pylyavskyy \cite{CP}.  We first give a 2-parking function of length $k$ from each 
3-spanning tree on $2k+1$ vertices arising from a fixed perfect matching.

\begin{figure}[ht]
\centerline{\psfig{file=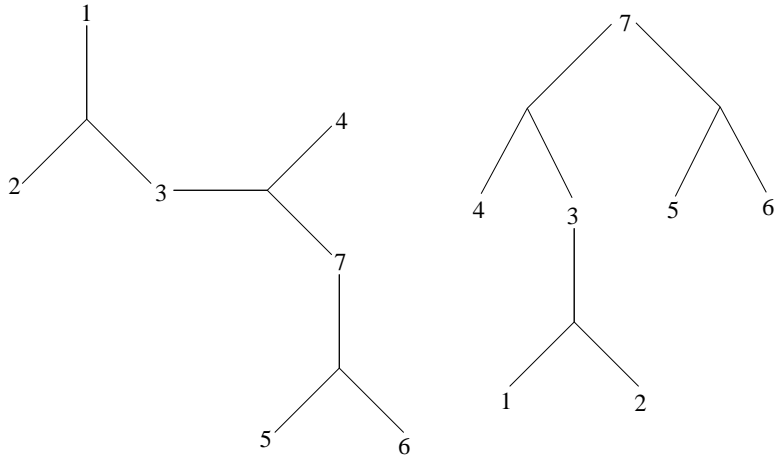,height=6cm}}
\caption{A spanning tree, redrawn with vertex 7 as the root.}
\label{fig:bfs}
\end{figure}

We use the 
BFS ordering of vertices of $[n]$.  We root the given 3-spanning tree at the vertex $n$ and 
write the other vertices in the order of their ``distance" from the root, breaking ties by
the natural order $1 < 2< \cdots < n$.  For example, given the spanning tree in Figure
\ref{fig:bfs}, we redraw it rooted at vertex 7 as in the digram on the right of 
Figure \ref{fig:bfs} and order the vertices as $7 < 3 < 4 < 5 < 6 < 1 < 2$.

\begin{proof}(Of Theorem \ref{thm:bij})  
We first deal with the case $r=2$.  
For $n=2k+1$, let $\cT_k$ be the set of 3-spanning trees arising from the perfect 
matching $M = \{\{1,2\},\{3,4\},\ldots,\{2k-1,2k\} \}$.  As there is a natural total order on the 
edges of this perfect matching, we label the edge $\{1,2\}$ as the first edge, $\{3,4\}$ as the 
second edge and so on.  Given a 3-spanning  
tree $T$, let $\pi_T$ be the BFS order on the vertices of $T$.  The map $\pi_T$ induces an order
on the set of hyperedges containing a perfectly matched edge.  For example, for the 3-spanning
tree of Figure \ref{fig:bfs}, the order on the hyperedges containing $\{3,4\}$ is 
$\{3,4,7\} <_{\pi_T} \{3,4,5 \} <_{\pi_T} \{3,4,6 \} <_{\pi_T} \{1,3,4 \} <_{\pi_T} \{2,3,4 \}$.
For the $i$-th edge of $M$, we refer to this order as ${\pi_T}(i)$.

We give a map $f_k: {\cT}_k \mapsto \Pk_k^2$ as follows.  Let $T$ be a 3-spanning tree arising
from $M$.  We give a number $a_i$ for the $i$-th edge of $M$ for each $1\leq i \leq k$.  
The sequence $\oa = (a_1,a_2,\ldots,a_k)$ is our candidate $f_k(T)$.  Let $a_i$ be the number 
of hyperedges which precede the unique hyperedge containing the $i$-th matched edge 
(ie $\{2i-1,2i\}$) in the order $\pi_T(i)$.

In the 3-spanning tree of Figure \ref{fig:bfs}, the sequence of numbers will be 
$(1,0,0)$ as $\{1,2,3 \}$ is the $2$-nd hyperedge in the order ${\pi_T}(1)$, 
$\{3,4,7 \}$ is the first edge in the 
order ${\pi_T}(2)$ and $\{5,6,7\}$ is the first edge in the order ${\pi_T}(3)$.

The map $f_k$ is clearly one-to-one.  We first show that the sequences obtained satisfy
the property that $f_k(i) \leq 2(i-1)$.  

Let $\ob$ be the weakly increasing rearrangement of the sequence $\oa$.  We recall that 
we obtained the sequence $\oa$ from $T \in \cT_k$.  We denote by ${\hgt}(i)$, 
the height of the $i$-th perfectly matched edge in the $n$-rooted version of $T$ 
(in our example, ${\hgt}(1) = 2, {\hgt}(2) = 1, {\hgt}(3) = 1$).  
Consider perfectly matched 
edges in increasing order of heights.  It is clear that there is at least one perfectly matched 
edge (say edge $i_1$ )whose height is 1.  It is easy to check that $a_{i_1} = 0$ and hence 
$b_1 = 0 < 1$.  The BFS order $\pi_T$ induces an order on the matched edges of the tree according
to the occurrence of the matched edge (or height).  For example, in the tree of Figure 
\ref{fig:bfs}, the order $\pi_T(M)$ on the matched edges is $\{3,4\} < \{5,6\} < \{1,2\}$.  
It is easy to see that the orders $\pi_T(M)$ and $\ob$ are identical.  Hence, the value
of the $i$-th element of $\ob$ is the hyperedge number of $i$-th matched edge in $\pi_T(M)$.
Clearly, this value is maximised when the tree has only one hyperedge at each height when this
value is at most $2(i-1)$.

We now show how to invert the map $f_k$. Let $\oa \in \Pk_k^2$.  
We identify the $i$-th coordinate with
the $i$-th edge of the perfect matching.  For $1 \leq j \leq k$, let $m(j)$ be the $j$-th 
perfectly matched edge.  We recall the total order on the edges of the perfect matching
$m(1) < m(2) < \cdots < m(k)$.  Let $\ob = (b_1 \leq b_2 \leq \cdots \leq b_k$ be the 
weakly increasing rearrangement of $\ob$ satisfying $b_i = b_j$ implies $m(i) < m(j)$.  We 
add hyperedges in increasing order of $\ob$ to construct a 3-uniform spanning tree rooted at
the vertex $n$.  Since $b_1 = 0$, there is atleast one perfectly matched edge $i_1$ such that 
$a_{i_1} = b_1$ and the first hyperedge we add has the $i_i$-th perfectly matched edge and $n$
(ie we add $\{2i_1-1,2i_1,n\}$).  If there are any more indices $j$ such that $b_j = 0$, we add
them too similarly.  Thus we can assume that we are at iteration $r$ with $b_r>0$.  
Since $b_r \leq 2(r-1)$, and since we have added $r-1$ perfectly matched pairs before the $r$-th
iteration, and since the vertex $n$ also exists, we have $2(r-1) + 1$ vertices already.  We can
consider the BFS order on the subtree consisting of just these $r-1$ hyperedges.  Let the
$b_r$-th vertex in the BFS be $v_{b_r}$.  Let $i_r$ be such that $a_{i_r} = b_r$.  We add the 
hyperedge $\{2i_r-1,2i_r,v_{b_r} \}$.  We do the same procedure till we get to index $p$ such
that $b_p > b_r$ (ie all of them are ``connected" to the same vertex in the tree).  Proceeding
this way we get a tree $T \in \cT_k$ as for each perfectly matched edge, we get a connection 
point which 
is already connected to the earlier tree.  Let $g_k(\oa)$ be the tree obtained by this procedure.
It is simple to check that $f_k(g_k(\oa) ) = \oa$.
\end{proof}

The proof for $r \geq 3$ is identical and following the proof above, we see that 
from $(r+1)$-spanning trees arising from a fixed $r$-perfect matching, we would 
get {\sl $r$-parking functions}.

\section{$r$-extended Shi hyperplane arrangements}
\label{sec:shi}
In this subsection, we point out that the results obtained above and a theorem of Stanley
\cite{stan-rota}, show that the number of $(r+1)$-spanning trees 
on $rk + 1$ vertices arising from a fixed $r$-perfect matching is equal to the number of 
regions of the $r$-extended Shi arrangement $S^r_k$.
We recall the following theorem from Stanley \cite{stan-rota}.


\begin{theorem}[Stanley\cite{stan-rota}]
\label{thm:stanley}
The number of $r$-parking functions of length $k$ is equal to the number of regions 
of the $r$-extended Shi arrangement $S_k^r$
\end{theorem}

The proof of Theorem \ref{thm:easy} is straightforward from the above theorem and Theorem
\ref{thm:bij}.
From \cite{stan-rota}, we see that the number of $r$-parking functions of length $k$ is
identical to the number of rooted $r$-forests on the vertex set $k$.  From Theorem \ref{thm:bij},
we see that both are identical to the number of $(r+1)$-spanning trees on $rk+1$ vertices
arising from a fixed $r$-perfect matching.  We note that
there is a simple bijection between such rooted $r$-forests on $[k]$ and 
$(r+1)$-spanning trees on $[rk+1]$ arising from a fixed $r$-perfect matching.

\section*{Acknowledgement}
We thank Professor Murali K. Srinivasan for illuminating discussions on the exponential
formula.

\bibliography{enumerate}
\end{document}